\documentstyle[12pt]{article}
\title{\bf Boundedness
and Stability of Impulsively Perturbed
Delay Differential Equations}

\author
{L. Berezansky $^{1}$
\\Ben-Gurion University of the Negev, \\
Department of Mathematics and Computer Science, \\
Beer-Sheva 84105, Israel, \\
E. Braverman $^{2}$ \\
Technion - Israel Institute of Technology,\\
Department of Mathematics, Haifa 32000, Israel }

\begin{document}
\maketitle

\footnotetext[1]{Supported by Israel Ministry of Science and
Technology}

\footnotetext[2]{Supported by the Centre for Absorption
in Science, Ministry of Immigrant Absorption State of Israel}

\section{Introduction}

{}~~~~~It is characteristic for a linear ordinary differential
equation that if any solution is bounded on the half-line
for any bounded right-hand side then
a solution of the corresponding
homogeneous equation tends to zero exponentially [1].
The connection of boundedness with exponential behavior of
solutions for impulsive differential equations is studied
in [2,3] and many other papers.
It turns out that for impulsive equations sometimes we
can avoid checking the boundedness of solutions for any
bounded right-hand side.
In particular the following
result is valid [4].

Suppose the points of impulses $\tau_j$ are such that
$0< \rho \leq \tau_{j+1}- \tau_j \leq \sigma < \infty $
and a solution of the scalar impulsive equation
\begin{eqnarray}
\dot{x}(t) + A(t)x(t)=0,
\end{eqnarray}
\begin{eqnarray}
x(\tau_j ) = B_j x(\tau_j -0) + \alpha_j, ~j=1,2, \dots ,
\end{eqnarray}
is bounded for each initial value $x(0)$
and each bounded sequence $\{ \alpha_j \}$.
(Here and in sequel we assume a solution is right continuous.)

Then the solution $X(t)$ of the homogeneous equation (1),
$x(\tau_j)=B_j x(\tau_j-0)$ satisfying
$x(0)=1$, has the exponential estimate
\begin{eqnarray}
\mid X(t) \mid \leq N e^{- \nu t},
\end{eqnarray}
with $N>0,~ \nu > 0$, and any solution of the non-homogeneous
equation
\begin{eqnarray}
\dot{x}(t) + A(t)x(t) = f(t),
\end{eqnarray}
is bounded for $t \geq 0$, if
$\sup_{t \geq 0} \mid f(t) \mid < \infty,
{}~ \sup_j \mid \alpha_j \mid < \infty . $

The first result of this type for the equation
\begin{eqnarray}
 \dot{x} (t) + A(t)x(t) = \sum_{j=1}^{\infty}
\alpha_j \delta(t-\tau_j) ,
\end{eqnarray}
where $\delta (t-\tau_j)$ are delta functions,
was obtained in [6].

This paper is concerned with the problem whether this result is valid
for impulsive delay differential equations.
The point is that the proof of this fact for (1),(2)
is essentially based on the representation [5]
\begin{eqnarray}
x(t) = X(t)x(0) +
\int_0^t G(t,s) f(s) ds + \sum_{0 < \tau_j \leq t}
G(t, \tau_j) \alpha_j
\end{eqnarray}
for any solution of (5),(2), where
\begin{eqnarray}
G(t,s)=X(t)X^{-1}(s).
\end{eqnarray}

For impulsive equations with delay we have a similar to (6)
solution representation ,
but (7) generally
speaking is not valid.
Nevertheless as we demonstrate
in this paper under certain conditions $\mid G(t,s)\mid
\geq \mid X(t)X^{-1}(s)\mid$.
Then the above result is valid for equations with delay.
It is to be emphasized that generally speaking the delay
is unbounded.

Besides this we obtain an explicit relation (Theorem 5)
connecting solutions of the impulsive equation with solutions
of the non-impulsive delay differential equation.
The corollary of Theorem 5 generalizes a result [7]
on the existence of a non-oscillating solution of an impulsive
delay differential equation.

\section{Preliminaries}

{}~~~~~Let
 ${\bf l}_{\infty}$ be a space of bounded sequences $\alpha = \{
 \alpha_j \}_{j=1}^{\infty} ,
 ~\alpha_j \in {\bf R},
{}~ \parallel \alpha \parallel_{{\bf l}_{\infty}} = \sup_j
\mid \alpha_j \mid $.

We consider a scalar linear delay differential equation
\begin{eqnarray}
\dot{x} (t) + \sum_{i=1}^m {A_i (t) x[h_i(t)]} =
r(t),~ t \geq 0,
{}~~~~x(\xi ) = \varphi (\xi),~ \xi < 0,
\end{eqnarray}
\begin{eqnarray}
x(\tau_j) = B_j x(\tau_j - 0) + \alpha_j, ~j=1,2, \dots ,
\end{eqnarray}
under the following assumptions

(a1) $ 0 = \tau_0 < \tau_1 < \tau_2 < \dots $
are fixed points,  $\lim_{j \rightarrow \infty} \tau_j = \infty $;

(a2) $A_i,r , ~i=1, \dots ,m $ are Lebesgue measurable
essentially bounded in any finite segment
$[0,b]$ functions,
$ ~B_j \in {\bf R}, j=1, \dots $ ;

(a3) $h_i: [0, \infty) \rightarrow {\bf R}$ are Lebesgue
measurable functions, $h_i(t) \leq t$ ;

(a4) $\varphi: (- \infty,0) \rightarrow {\bf R}$
is a Borel measurable bounded function.

Below we will also need the following hypotheses

(a5) $M = \sup \parallel B_j \parallel < \infty $;

(a6) there exists $\rho > 0$ such that $\tau_{j+1} - \tau_j \geq
\rho , ~j=1,2, \dots $ ;

(a7) there exists $\sigma > 0$
such that $\tau_{j+1} - \tau_j \leq \sigma, ~j=1,2, \dots $;

(a8) there exists $\delta > 0$ such that
$t - \delta \leq h_i (t) , ~i=1, \dots, m$.

(a9) there exists $Q>0$ such that
$$ \int_k^{k+1} \mid A_i (t)\mid ~dt \leq Q ,
{}~k = 1,2, \dots ,~i=1, \dots, m . $$

\underline{\sl Definition.}
A function $x: [0, \infty) \rightarrow {\bf R}$ absolutely
continuous in each $[\tau_j, \tau_{j+1} ) $
is {\bf a solution} of the impulsive equation (8),(9),
if for $t\neq\tau_j$ it satisfies (8) and for $t=\tau_j$
it satisfies (9).

A solution $X(t)$ of the homogeneous equation
\begin{eqnarray}
\dot{x} (t) + \sum_{i=1}^m {A_i (t) x[h_i(t)]} =
0
\end{eqnarray}
for $t \geq 0$,
\begin{eqnarray}
x(\xi) = 0
\end{eqnarray}
for $\xi < 0$, satisfying impulsive conditions
\begin{eqnarray}
x(\tau _j) = B_j x(\tau _j - 0),
\end{eqnarray}
$j=1,2, \dots ,$
and such that $x(0)=1$ is said to be
{\bf a fundamental solution}.

A solution $G(t,s)$ of the homogeneous equation
(10), $ t \geq s, $
(11), $  \xi < s,$ (12), $\tau_j > s$,
satisfying $G(s,s) = 1$
is said to be {\bf a fundamental function }
(sometimes it is called a Green or an evolutionary function).
We assume $G(t,s) = 0, ~0 \leq t < s$.

\newtheorem{guess}{Theorem}
\begin{guess}
{\bf [8]} Suppose (a1)-(a4) are satisfied.
Then there exists one and only one
solution of the problem (8) , with the initial
value
$$x(0) = \alpha _0
$$
and impulsive conditions (9)
and it can be presented as
\begin{eqnarray}
x(t) = X(t)x(0) + \int_0^t {G(t,s) r(s) ds} - \nonumber \\
- \sum_{i=1}^m {\int_0^t {G(t,s) A_i(s) \varphi (h_i(s)) ds}} +
\sum_{ 0 < \tau_j \leq t} {G(t, \tau_j) \alpha _j} .
\end{eqnarray}
Here $\varphi (\zeta) = 0,$ if $\zeta \geq 0, $
and $X(t) = G(t,0)$.
\end{guess}

The fundamental function $G(t,s)$ generally speaking does not satisfy (7).
However in certain cases the inequality
$ \mid G(t,s) \mid
 \geq \mid  X(t) X^{-1} (s) \mid $
holds as in the case of delay differential equations
without impulses [9].

\newtheorem{uess}{Lemma}
\begin{uess}
Suppose (a1)-(a4) hold.
Let  $G(t,s) \geq 0, ~X(t) > 0,~
A_i(t) \geq 0 $ for any $t,s, ~ 0 \leq s \leq t < \infty,
{}~i=1, \dots ,m. $

Then $G(t,s) \leq G(t,\zeta) G(\zeta, s)$ and $G(t,s) \geq X(t)
X^{-1} (s),~0\leq s\leq\zeta <t$.
\end{uess}

{\sl Proof.}
Let $0 < s < \zeta$ .
Suppose $y$ is a function satisfying (10),(12)
for $t> \zeta, ~\tau_j > \zeta$.
Besides this, let $y(\xi) = G(\xi,s)$,
if $s< \xi < \zeta, ~ y(\xi)=0$, if $\xi < s$, and
$y(\zeta) = G(\zeta,s) .$
Then by Theorem 1
$$ G(t,s) = y(t) = G(t, \zeta) G(\zeta, s) -
\sum _{i=1}^m \int_{\zeta}^t G(t, \xi) A_i (\xi)
G(h_i (\xi), s) d \xi ,$$
where $G(h_k(\xi),s) = 0$,
if $h_k(\xi)<s$ or $h_k(\xi) > \zeta$.
As $G(t,s) \geq 0, ~ A_i (t) \geq 0 ,$ then the second
term in the right-hand side is not positive.
Therefore
$G(t,s) \leq G(t,\zeta) G(\zeta, s),
s \leq \zeta \leq t, $
$$ X(t)X^{-1}(s) = G(t,0)G^{-1}(s,0) \leq
G(t,s)G(s,0)G^{-1}(s,0) = G(t,s),$$
which completes the proof.

\section{Main results}

\begin{guess}
Suppose $A_i (t) \geq 0, ~G(t,s) \geq 0, ~X(t) > 0,$
the hypotheses (a1)-(a4),(a7) are satisfied,
and the solution of the problem (10),
$t \geq 0,$ (11), $ \xi < 0, $ (9), $x(0) = 0$
is bounded for each
$\alpha = \{ \alpha_i \}_{i=1}^{\infty} \in {\bf l}_{\infty}.$

Then there exist positive constants $N$ and $\nu$
such that (3) holds.
\end{guess}

{\sl Proof.}
By Theorem 1 the solution of the problem (10),
$t \geq 0$, (11), $\xi < 0$, (9),
$x(0)= 0$ can be presented as
$$ x(t) = \sum_{0 < \tau_i \leq t} G(t, \tau_i ) \alpha_i .$$
The right hand side of this formula for any $t$
defines a bounded linear operator acting from ${\bf l}_{\infty}$
to ${\bf R}$ since
$$ \mid x(t) \mid \leq
\sum_{0 < \tau_i \leq t} \mid  G(t, \tau_i) \mid
\parallel  \alpha \parallel _{{\bf l}_{\infty}} , $$
where $\alpha = (\alpha_1, \alpha_2, \dots , \alpha_i , \dots )
\in {\bf l}_{\infty} $.
By the hypothesis of the theorem for any sequence
$\alpha \in {\bf l}_{\infty}$ the solution $x(t)$ is bounded.
Therefore
the uniform boundedness principle implies
that there exists $k > 0$ such that
$$ \mid \sum_{0 < \tau_i \leq t} G(t, \tau_i ) \alpha_i \mid
\leq  k \parallel \alpha \parallel _{ {\bf l}_{\infty}} .
$$
Without loss of generality we can assume $k>1$.
Since $G(t,s) \geq 0 $ then
$$  \sum_{0 < \tau_i \leq t} G(t, \tau_i ) \leq k  .$$
By Lemma 1
$~~ G(t, \tau_i) \geq X(t)X^{-1} (\tau_i) .~ $
Thus for $t= \tau_2$
$$G(\tau_2, \tau_2) + G(\tau_2, \tau_1) =
1 + G(\tau_2 , \tau_1) \leq k $$
implies
$$ \frac{X(\tau_2)}{X(\tau_1)} \leq k-1 , $$
therefore
$ X(\tau_2) \leq (k-1) X(\tau_1) . $
By assuming
\begin{eqnarray}
X(\tau_i) \leq X(\tau_1) (k-1)^{i-1} /k^{i-2},
\end{eqnarray}
$2 \leq i \leq j, $ we obtain
$$ k \geq \sum_{i=1}^{j+1} G(\tau_{j+1} , \tau_i)
\geq 1 + \sum_{i=0}^{j-1} \frac{X(\tau_{j+1})}{X(\tau_{i+1})} \geq
$$ $$\geq 1 +
 \frac{X(\tau_{j+1})}{X(\tau_1)} \left[ 1 +
\sum_{i=0}^{j-2}
\frac{k^i}{(k-1)^{i+1}} \right] =
 1 +
 \frac{X(\tau_{j+1})}{X(\tau_1)} \left[ 1 +
\frac{ \frac{k^{j-1}}{(k-1)^j} - \frac{1}{k-1} }{\frac{k}{k-1} - 1}
\right] = $$
$$  = 1 +
 \frac{X(\tau_{j+1})}{X(\tau_1)} \left[ 1 +
\frac{k^{j-1}}{(k-1)^{j-1}} - 1 \right] =
1 +
 \frac{X(\tau_{j+1})}{X(\tau_1)} \frac{k^{j-1}}{(k-1)^{j-1}} . $$
Hence
$$
 \frac{X(\tau_{j+1})}{X(\tau_1)} \leq
\frac{(k-1)(k-1)^{j-1}}{k^{j-1}} =
\frac{(k-1)^j}{k^{j-1}}. $$
Thus by induction we obtain (14) for any $i = 3,4, \dots . $

Let $\tau_j < t \leq \tau_{j+1}$.
As $X(t) > 0, ~A_i (t) \geq 0,$
then in $[\tau_j, \tau_{j+1} )$ $X(t)$ does not increase.
Since by (a7)
$~~t \leq (j+1) \sigma$, i.e. $j \geq t/\sigma-1$, then
$$\ln X(t) \leq \ln X(\tau_j)
\leq \ln X(\tau_1) - (t/\sigma - 2) \ln [k/(k-1)]
+  \ln k. $$
Therefore (3) holds, with
$$ \nu = \ln [k/(k-1)] / \sigma , ~~
 N = \max \left\{ X(\tau_1) k^3 / (k-1)^2 , ~
\sup_{0 \leq t \leq \tau_1} [\exp (\nu t ) X(t) ] \right\}. $$

\begin{guess}
Suppose the hypotheses of Theorem 2 and
(a5),(a9) hold.

Then there exist positive constants $N$ and $\nu$ such
that
\begin{eqnarray}
\mid G(t,s) \mid \leq N \exp [-\nu (t-s)].
\end{eqnarray}
\end{guess}

{\sl Proof.}
Similar to the proof of Theorem 2 we obtain
$$0 \leq G(t,s) \leq N_s \exp [- \nu (t-s)],
\mbox{~~with~~}
\nu = \ln [k/(k-1)] / \sigma , $$
$$ N_s = \max \{ G(\tau_p ,s) k^3 / (k-1)^2 , ~
\sup_{s \leq t \leq \tau_p} [\exp (\nu (t-s) ) G(t,s) ] \} , $$
where $\tau_p$ is the least of $\tau_j > s.$
Hence $\tau_p - s \leq \sigma $ and
$$\sup_{s \leq t \leq \tau_p} [\exp (\nu (t-s) ) G(t,s) ] \leq
\exp (\nu \sigma )
\sup_{s \leq t \leq \tau_p} G(t,s)   . $$
In Lemma 3.2 of the paper [8] the following estimate is
established
$$ \mid G(t,s) \mid \leq
(1 + \mid B_p \mid )
\exp \left\{ \int_s^t \sum_{i=1}^m
\mid A_i (\zeta) \mid d \zeta \right\} , $$
$\tau_{p-1} < s \leq t \leq \tau_p$ .
Therefore
$$ G(t,s) \leq (1+M) \exp [mQ (t-s)] \leq (1 + M) \exp (mQ \sigma ) $$
for $\tau_{p-1} < s \leq t \leq \tau_p .$
Hence by assuming
$$ N = \max \left\{
\frac{(1+M) \exp(mQ \sigma) k^3}{(k-1)^2 } ,
(1+M) \exp [ \sigma (\nu + mQ) ] \right\}, $$
we obtain the estimate (15) for $G(t,s)$.

\begin{guess}
Suppose (a1)-(a9) are satisfied, $A_i (t) \geq 0, ~G(t,s) \geq 0,
{}~X(t) > 0$
and the solution of (10), $t \geq 0$,
(11),$\xi > 0$, (9),  $x(0)=0$ is bounded on $[0, \infty)$
for any $\alpha = \{ \alpha_i \}_{i=1}^{\infty} \in
{\bf l}_{\infty}$.

(a)  If $\sup \mid \alpha_i \mid < \infty,
{}~ \sup_{t \geq 0} \mid r(t) \mid < \infty,$
then any solution of (8),(9) is bounded on the half-line $[0,
\infty)$ ;

(b)  if $\lim_{n \rightarrow \infty} \mid \alpha_n \mid = 0,
{}~ \lim_{t \rightarrow \infty} \mid r(t) \mid  = 0$,
then for any solution $x$ of (8),(9) $\lim_{t \rightarrow \infty}
\mid x(t) \mid = 0$;

(c)  if there exist positive constants $P$ and $\lambda$
such that $\mid \alpha_n \mid \leq P e^{- \lambda n}$ and $
\mid r(t) \mid ~ \leq  P e^{- \lambda t}, $ then for any
solution $x$ of (8),(9) there exist $P_0>0,~\lambda_0 >0$
such that $\mid x(t) \mid \leq P_0 e^{- \lambda_0 t} . $
\end{guess}

{\sl Proof.}
The similar result for impulsive equations without delay is
obtained in [4].
Comparing solution representations (6) and (13) we obtain that
only the terms
$$ \int_0^t G(t,s) A_i (s) \varphi (h_i (s)) ds,
{}~ \varphi(\xi) = 0, \xi \geq 0, $$
have to be estimated.
Since by (a8) $\varphi(h_i(s)) = 0, s \geq \delta $,
 then by  Theorem 3
$$ \mid \int_0^t G(t,s) A_i (s) \varphi(h_i(s)) ds
\mid
\leq \int_0^{\delta} \mid G(t,s) \mid
\mid  A_i(s) \mid  \mid  \varphi (h_i(s))
\mid  ~ds \leq $$
$$ \leq \sup_{s \in [0, \delta ] } \!\! \mid \!\!
A_i (s) \! \mid \! \sup_{s<0} \!
\mid \! \! \varphi (s) \! \! \mid \!\int_0^{\delta}
\!\! N e^{- \nu (t-s)} ds \leq
\sup_{s \in [0, \delta ]} \!\! \mid \!
\! A_i (s) \! \! \mid  \sup_{s<0}
\mid \!\! \varphi (s) \!\! \mid \!
\frac{N}{\nu} e^{ \nu \delta} e^{- \nu t} ,$$
which is obviously bounded, tends to zero
as $t$ tends to infinity and has an exponential estimate.

\underline{\sl Remark.}
(c) yields that under the hypotheses of Theorem 4 the equation
(10),(12) is exponentially stable, i.e. for its solution $x$
the following inequality holds
$$\mid\! x(t) \!\mid \leq N e^{- \nu t}
\left( \mid\! x(0) \!\mid + \sup_{s<0}  \!\mid \varphi(s) \!\mid \right) . $$
\vspace{5 mm}

All the hypotheses of Theorems 2-4
are easily verified except the conditions
$ G(t,s) \geq 0 , ~X(t) > 0. $
For $G(t,s)$ the representation (7) is not valid.
However sufficient conditions for $C(t,s)$ being positive
are known [10,11], where $C(t,s)$ is a fundamental function
of the equation (8) without impulses.
For instance, $C(t,s)$ is positive if
\begin{eqnarray}
\sum_{i=1}^m \int_{h_i^+(t)}^t A_i^+ (s) ds \leq 1/e,
\end{eqnarray}
where $a^+ = \max \{ a,0 \} $ .

We study the following problem:
is $G(t,s)$ positive for
$C(t,s)$ being positive ?
Generally speaking it is not true.
For an ordinary differential equation this assertion is valid
if and only if $B_j > 0, ~j = 1,2, \dots .$
For an impulsive delay differential equation
the positiveness of $C(t,s), ~B_j, j = 1,2, \dots , $
does not imply the positiveness of $G(t,s)$.

{\bf Example.}
By (16) the equation
$$\dot{x} (t) + x(t- \frac{1}{3} ) = f $$
has $C(t,s) > 0$ .
Consider this equation with impulsive conditions
$$x(\frac{j}{3}) = \frac{1}{6} x(\frac{j}{3}-0), ~j=1,2, \dots . $$
Then
$G(t,0) = 1, ~t \in [0, \frac{1}{3} ), $
$ G(t,0) = \frac{1}{6} - (t- \frac{1}{3}), ~t\in [\frac{1}{3} ,
\frac{2}{3} ). $
Thus
$G(t,0) < 0$ for $ \frac{1}{2} < t < \frac{2}{3} .$

\vspace{5 mm}

In many cases we can construct the fundamental function $C(t,s)$ of the
non-impulsive equation (8).
Now we will obtain a relation connecting the fundamental function
of the impulsive equation $G(t,s)$, with $C(t,s)$.
The corresponding relation for impulsive equations without
delay is [5]
\begin{eqnarray}
G(t,s) = \left\{
\begin{array}{l}
C(t,s),~\tau_i \leq t,s < \tau_{i+1}, \\
C(t, \tau_i) \left[
\prod_{j=k+1}^{i} B_j C(\tau_j, \tau_{j-1}) \right]
B_k  C(\tau_k, s) , \\
\tau_{k-1} \leq s < \tau_k \leq \tau_i \leq t < \tau_{i+1}.
\end{array} \right.
\end{eqnarray}

Let $k$ and $j$ be positive integers, $k \leq j$.
Denote by $\Omega_{k,j}$ the set of all non-empty ordered subsets
$e= \{ n_1, \dots, n_i \}$ of the set
$\{ k, k+1, \dots, j-1, j \}$, with natural order.
Denote by $\max(e)=n_i$ and $\min(e)=n_1$
a minimal and a maximal elements of $e$ correspondingly.
Let for any $e \in \Omega_{k,j}$  a set $\Lambda_e$ be
the corresponding set of pairs
$\{ (n_i,n_{i-1}), \dots , (n_2,n_1) \}$,
$\Lambda_e = \emptyset $ if $e$ is a one-point set.

\begin{guess}
Suppose (a1)-(a4) hold and
$C(t,s)$ is a fundamental function of the non-impulsive
delay differential equation (8).

Then
$~G(t,s) = C(t,s), ~~\tau_l \leq t,s < \tau_{l+1}$,
\begin{eqnarray}
G(t,s) = C(t,s) + \!\!\! \sum_{e \in \Omega_{k,l}}
\!\!\!\! C(t, \tau_{\max(e)}) \left[  \prod _{(n_p,n_{p-1}) \in \Lambda_e}
\!\!\!\!\!\!\!\! (B_{n_p} - 1)
 C(\tau_{n_p}, \tau _{n_{p-1}}) \right] \times \nonumber \\
\times  (B _{\min (e)} - 1) C(\tau_{\min(e)},s) ,
{}~~~~\tau_{k-1} \leq s < \tau_k \leq \tau_l \leq t <
\tau_{l+1}.~~~~~~
\end{eqnarray}
Here $\prod _{(n_p,n_{p-1}) \in \Lambda_e} = 1$ if $e$ is a
one-point set.
\end{guess}

{\sl Proof.}
For $s < t < \tau _1$ solution $G(t,s)$ of the impulsive equation
and $C(t,s)$ of the equation without impulses obviously coincide.
If $t=\tau _1$, then
$$ G(\tau _1,s) = B_1 C(\tau _1,s).$$
Let $0 \leq s <\tau _1 < t < \tau _2.$
Then the solution $G(t,s)$ of the impulsive equation (10),
$t \geq s$, (12)
can be treated as the solution of the equation (10) on
the segment $[\tau _1, \tau _2]$ with the initial value
$ G(\tau _1,s) = B_1 C(\tau _1, s)$
and the initial function
$C(t,s).$
By Theorem 1
\begin{eqnarray}
G(t,s) = C(t,\tau_1) B_1 C(\tau_1,s) -
\sum _{i=1}^m \int_{\tau _1}^t C(t,\xi) A_i (\xi) C(h_i(\xi),s)
{}~d\xi ,
\end{eqnarray}
where $C(h_i(\xi),s)=0,$~ if $ h_i (\xi) > \tau_1 $.

On the other hand, $C(t,s)$ for $0 \leq s <\tau_1 < t < \tau_2$
can also be treated as a solution of (10) on the segment $[\tau_1,\tau_2]$
with the same initial function $C(\zeta ,s), \zeta \in [s, \tau_1 ) $
and the initial value $C(\tau_1,s).$
Therefore
\begin{eqnarray}
C(t,s) = C(t,\tau_1) C(\tau_1,s) -
\sum _{i=1}^m \int_{\tau _1}^t C(t,\xi) A_i (\xi) C(h_i(\xi),s) d
\xi.
\end{eqnarray}

By subtracting (19) from (18) one obtains
$$ G(t,s) - C(t,s) = C(t,\tau_1)(B_1 - 1) C(\tau_1,s),$$
hence for $0 \leq s < \tau_1 \leq t < \tau_2$
$$ G(t,s) = C(t,s) + C(t,\tau_1)(B_1 - 1) C(\tau_1,s).$$

In fact the latter formula defines how the fundamental function
changes when it overcomes the point $\tau_1$.
Let $G_l (t,s)$ be a function
satisfying the homogeneous equation (10), $t \geq s$,  and the impulsive
conditions
$$G_l(\tau_i, s) = B_i G_l(\tau _i -0, s),~ i = 1, \dots, l-1,
{}~\tau_i > s.$$
Obviously $G_l(t,s) = G(t,s) ,~ t < \tau _l.$
Then for $0 \leq s < \tau_l \leq t < \tau _{l+1}$ we similarly obtain
\begin{eqnarray}
G(t,s) = G_l(t,s) + G_l(t,\tau_l)(B_l - 1) G_l(\tau_l,s) .
\end{eqnarray}

We prove the equality (18) by induction in $l$.
For $l=1$ this equality was obtained above.

Suppose that for $0 \leq \tau_{k-1} \leq s
< \tau_k \leq \tau_{l-1} \leq t < \tau _l$ this equality is valid.
Then $G_l (t,s)$
for $t \in [\tau_l, \tau _{l+1} )$ can be presented as
$G_l (t,s) =  C(t,s) + $ $$ + \!\!\! \sum_{e \in \Omega_{k, l-1}}
\!\!\!\!\! C(t, \tau_{\max(e)})
 \left[  \prod _{(n_p,n_{p-1}) \in \Lambda_e}
\!\!\!\!\!\!\! (B_{n_p} - 1) C(\tau_{n_p}, \tau _{n_{p-1}} ) \right]
(B_{\min(e)} - 1) C(\tau_{\min(e)} ,s) . $$
Let $0 \leq \tau_{k-1} \leq s < \tau_k \leq \tau_l \leq t < \tau_{l+1}$.

Obviously $\Omega_{k,l} = \Omega_{k,l-1} \cup \{ l \}
\cup \left\{  \{ e , l \}  \mid  e \in \Omega_{k,l-1}
\right\} $ and $G_l(t, \tau_l) = C(t, \tau_l)$.

Therefore by substituting this in (21) we obtain
$$G(t,s) = G_l(t,s) + G_l(t,\tau_l)(B_l - 1) G_l(\tau_l,s)
=
C(t,s) + $$ $$ + \!\!\! \sum_{e \in \Omega_{k,l-1}}
\!\!\!\!\! C(t, \tau_{\max(e)}) \left[ \! \prod _{(n_p, n_{p-1}) \in \Lambda_e}
\!\!\!\!\!\!\!\! (B_{n_p} - 1) C(\tau_{n_p}, \tau _{n_{p-1}} ) \right]
\! (B_{\min(e)} - 1) C(\tau_{\min(e)},s)  + $$
$$ + C(t, \tau _l)(B_l - 1)C(\tau_l,s)
+ \sum_{e \in \Omega_{k,l-1}}
\!\! C(t, \tau_l)(B_l - 1)
C(\tau_l , \tau_{\max(e)}) \times $$
$$ \times \left[  \prod _{(n_p, n_{p-1}) \in \Lambda_e}
\!\!\! (B_{n_p} - 1) C(\tau_{n_p}, \tau _{n_{p-1}} ) \right]
(B_{\min(e)} - 1) C(\tau_{\min(e)},s)
= C(t,s) + $$ $$  + \!\! \sum_{e \in \Omega_{k,l}}
\!\! C(t, \tau_{\max(e)}) \left[  \prod _{(n_p, n_{p-1}) \in \Lambda_e}
\!\!\!\!\!\! (B_{n_p} - 1) C(\tau_{n_p}, \tau _{n_{p-1}} ) \right]
(B_{\min(e)} - 1) C(\tau_{\min(e)},s)  . $$

The proof of the theorem is complete.

\underline{\sl Corollary 1.}
Suppose the fundamental function of the non-impulsive equation
$C(t,s) > 0 ,~ B_i \geq 1, ~ i =1,2, \dots . $

Then the fundamental function of the impulsive equation
$G(t,s)>0$.

\underline{\sl Corollary 2.}
Suppose (16) holds, $B_j \geq 1$.
Then the fundamental solution is non-oscillating.

{\bf Remarks.}
1. By the straightforward verification one can prove
that (18) implies (17) in the case of an impulsive
differential equation without delay.

2. In the vector case Theorem 5 is valid if in (18) we change
unit by the unit matrix.

3. Corollary 2 generalizes the result obtained in [7]
on the existence of a non-oscillating solution of an impulsive
equation with constant delay and constant coefficients.

\end{document}